\documentclass[12pt]{amsart}
\usepackage{amsfonts,amsmath,amssymb,amscd}
\newtheorem{theorem}{Theorem}
 
\newtheorem{lemma}[theorem]{Lemma} 
 
\newtheorem*{cor}{Corollary} 
\newtheorem{prop}[theorem]{Proposition} 
 
\newtheorem*{main}{Main Theorem}
\newtheorem*{S-vN}{Stone-von Neumann Theorem}
\DeclareMathOperator{\GL}{GL}
\DeclareMathOperator{\PGL}{PGL}
\DeclareMathOperator{\Sp}{Sp}
\DeclareMathOperator{\Gal}{Gal}
\DeclareMathOperator{\Aut}{Aut}
\DeclareMathOperator{\End}{End}
\DeclareMathOperator{\cha}{char}
\DeclareMathOperator{\im}{im}
\DeclareMathOperator{\supp}{supp}

\renewcommand{\O}{\mathcal O}
\newcommand\cO{{\mathcal O}} 
\newcommand{\maxi}{\mathfrak m}
\newcommand{\id}{\mathfrak i}
\newcommand{\pd}{\mathfrak p}

\newcommand\idl{\id_\lambda} 
\newcommand\rn{{\bf Q}} 
\newcommand\Q{{\bf Q}} 
\newcommand\Z{{\bf Z}}
\newcommand\C{{\bf C}}  
\newcommand\cn{\bf C} 
\newcommand\cS{\mathcal S} 
\newcommand\cV{\mathcal V} 
\newcommand\aform[2]{\left<#1,#2\right>}

\begin{document}

\title{Realizing the local Weil representation over a number field}
\author{Gerald Cliff and David McNeilly}
\address{University of Alberta, Department of Mathematical and Statistical
Sciences, Edmonton, Alberta, Canada T6G 2G1}
\email{gcliff@math.ualberta.ca, dam@math.ualberta.ca}
\begin{abstract}We show that the Weil representation of the symplectic
group $\Sp(2n,F)$, where $F$ is a non-archimedian local field, can be
realized over the field $K=\Q(\sqrt p,\sqrt{-p})$, where $p$ is the
residue characteristic of $F$. 
\end{abstract}

\maketitle

\section{Introduction}
Our main result is that the Weil representation of the symplectic
group $\Sp(2n,F)$, where $F$ is a non-archimedian local field of residue
characteristic~$\ne 2$, can be realized over a number field $K$. 
We take
an infinite-dimensional complex
vector space $\cV$ such that the  Weil representation is given by 
$\rho:\Sp(2n,F)\to\PGL(\cV)$
and we find a $K$-subspace $\cV_0$ of $\cV$ such that $\rho(g)(\cV_0)=\cV_0$
for all $g\in \Sp(2n,F)$.

This answers a question raised
by D.~Prasad [P]. Indeed, we show that we can take $K=\Q(\sqrt p,\sqrt{-p})$
where $p$ is the residue characteristic of $F$. We assume that $p$ is odd.
A consequence of this, also pointed out by Prasad, is that the local
theta correspondence can be defined for representations which are realized
over $K$.

Let $W$ be the Weil representation of $\Sp(2n,F)$.
The Weil representation can be defined using the Schr\"odinger
representation of the Heisenberg group $H$.
Let $\lambda$ be a fixed complex character on the additive group of
the field $F$. Suppose that $F^{2n}$ is the direct sum $X\oplus Y$ of
totally isotropic $F$-subspaces.
The Schr\"odinger model is realized in the
Bruhat-Schwartz space $\cS(X)$ of locally constant functions $f:X\to \C$
of compact support. For $h\in H$, there are operators 
$S_\lambda(h)$ on $\cS(X)$ such that $S_\lambda:H\to \GL(\cS(X))$ is the
unique smooth irreducible representation of $H$ with
central character $\lambda$. The natural action of the symplectic groups
extends to an action on~$H$, and the Weil representation is given
by operators $W_\lambda(g)$ on $\cS(X)$, $g\in \Sp(2n,F)$, such that
$$W_\lambda(g)^{-1}S_\lambda(h)W_\lambda(g)=S_\lambda(hg)\quad
h\in H,\,g\in \Sp(2n, F).$$

Let $\Q(\lambda)$ be the field obtained by adjoining
all the character values of $\lambda$ to $\Q$, and let
$E=\Q(\lambda)(\sqrt{-1})$. In the case that
$F$ has characteristic 0, $E$ is the field obtained from $\Q$ by
adjoining $\sqrt{-1}$  and all $p$-power roots of unity.
For a subfield $L$ of $\C$,
define $\cS(X,L)$ to be the space of locally constant 
functions on $X$ of compact support having values in $L$. 
We show that there is an explicit choice of Weil 
operators $W_\lambda(g)$ on $\cS(X)$ which leave  $\cS(X,E)$ invariant.


The Galois group of $E$ over $\Q$ acts on $\cS(X,E)$ and on
$\End(\cS(X,E))$.  In Section 7 we define a 1-cocycle $\delta$ on
$\Gal\bigl(E/\Q(\sqrt p,\sqrt{-p})\bigr)$ with values in
$\GL\bigl(\cS(X,E)\bigr)$ such that
\begin{equation}\label{delta0}
{}^\sigma W_\lambda(g)=\delta(\sigma)^{-1}W_\lambda(g)\delta(\sigma),\quad
g\in\Sp(V).
\tag{I}
\end{equation}
Using Galois descent, we show that there exists $\alpha\in\GL(\cS(X,E))$
such that $\delta(\sigma)=\alpha^{-1}\,{}^\sigma\!\alpha$ for $\sigma\in
\Gal\bigl(E/\Q(\sqrt p,\sqrt{-p})\bigr)$. 
\begin{main}The operators 
$\alpha W_\lambda(g)\alpha^{-1}$ leave $\cS(X,\Q\bigl(\sqrt p,\sqrt{-p})\bigr)$
invariant, and provide a form of the Weil representation realized
over $\Q(\sqrt p,\sqrt{-p})$.
\end{main}

To indicate how we find the 1-cocycle sastifying (I), for the rest of the
introduction we assume that $F$ has characteristic 0. The Galois
group of $\Q(\lambda)/\Q$ is isomorphic to the units $\Z_p^*$ of the
$p$-adic integers.  For an element $s$ of $\Z_p^*$, we let $\sigma_s$
denote the corresponding element of $\Gal(\Q(\lambda)/\Q)$.  For an
element $t\in F^*$, we define the character $\lambda[t]$ of $F$ by
$\lambda[t](r)=\lambda(tr)$, $r\in F$. 

For~$t \in F^*$, let~$g_t \in \Sp(2n,F)$ be defined by~$(x+y)g_t = t^{-1}x + ty$,
and~$f_t \in \GL(2n,F)$ by~$(x+y)f_t = x + ty$, where~$x \in X$,~$y \in Y$.
Then~$f_t$ is not in general in~$\Sp(2n,F)$, but conjugation by~$f_t$ leaves~$\Sp(2n,F)$ invariant. We have
\begin{equation}
W_\lambda(g^{f_t})=W_{[\lambda(t)]}(g),\quad g\in \Sp(V).
\tag{II}
\end{equation}
Furthermore, observing~$f_{t^2}$ is the composite~$tI \circ g_t$, we show
\begin{equation}
W_\lambda(g^{f_{t^2}})= W_\lambda(g_t)^{-1} W_\lambda(g) W_\lambda(g_t).
\tag{III}           
\end{equation}

For $\sigma\in\Gal(E/\Q(\sqrt p,\sqrt{-p}))$, $\sigma|_{\Q(\lambda)}$ 
is a square, so $\sigma|_{\Q(\lambda)}=\sigma_{\epsilon^{2i}s}$ where 
$\epsilon$ in a primitive $p-1$ root of unity, $s$ is a principal unit of $\cO$,
and $i$ is an integer, $1\le i\le(p-1)/2$.  We note
\begin{equation}
{}^\sigma W_\lambda(g)=W_{\lambda[\epsilon^{2i}s^2]}(g).
\tag{IV}
\end{equation}
In light of~(II) and~(III), we deduce
$$^{\sigma}W_\lambda(g) =
             W_\lambda(g_{\epsilon^is})^{-1} 
             W_\lambda(g) W_\lambda(g_{\epsilon^is}). 
$$     
The last equation is used to show that ~$\delta(\sigma) \
= W_\lambda(g_{\epsilon^is})$ satisfies (\ref{delta0}) and almost satisfies 
the one-cocycle condition.

Equations (II),~(III), and~(IV) are proved using an integral formula for
Weil operators due to Ranga Rao [RR]; see equation (3) in Section 3. 
This formula is also used
to show that conjugation by the Weil operators $W_\lambda(g)$ leaves
$\cS(X,E)$ invariant.

\section{Preliminary remarks on local fields, characters and measures}

We fix some notation and recall some elementary facts about the 
characters of the additive group of a local field. Further details
can be found in the first two chapters of~[W].

Let $F$ be a non-Archimedean local field, $\O$ its ring of integers, 
and $\maxi$ the maximal ideal of $\O$. The order of the residue class 
field $\kappa = \O/\maxi$ shall be denoted $q$; we note that $q$ is 
power of $p = \cha \kappa$. We assume throughout that $p$ is different 
from $2$; in particular, $2$ is a unit of $\O$. 

Given a fractional $\O$-ideal ${\mathfrak a}$, there
exists an unique integer $v({\mathfrak a})$, the~{\it valuation of $\mathfrak a$}, 
such that
            $$ {\mathfrak a} = \maxi^{v({\mathfrak a})}.$$
If $s \in F$ is non-zero, the valuation of the ideal $s\O$
is refered to as the valuation of $s$, denoted $v(s)$. 
The absolute value on $F$ is  related to the valuation $v$ on $F$ by
$$ \vert s \vert = q^{-v(s)},\quad s\in F, s\ne 0. $$

Let $\lambda$ be a non-trivial, continuous, complex linear (unitary)
character of $F^+$.  The continuity of $\lambda$ ensures 
that its kernel contains a fractional $\O$-ideal. The fact that $\lambda$
is non-trivial allows one to deduce that the set of all such
fractional $\O$-ideals has a unique maximal element $\id =
\id_\lambda$, {\it the conductor of $\lambda$}. The {\it level
of $\lambda$} is defined to be the valuation of $\id_\lambda$.
\par

Given $n \ge 1$, let
         $$ \nu_{p^n} = \{\, z \in {\bf C} \, : \, z^{p^n} = 1\, \},\quad
\nu_{p^\infty}=\bigcup_{n=1}^\infty\nu_{p^n}.$$
(The more customary symbol $\mu$ will be used to denote a measure.)
\begin{lemma}We have
\begin{equation*} \im \lambda = 
\begin{cases} \nu_p, & \text{if $\cha F =  p$;}\\
               \nu_{p^\infty}, & \text{if $\cha F = 0$.}
\end{cases}\end{equation*}
\end{lemma}
\begin{proof}
Take $x \in F$. If $\cha F =
p$ then
         $$ 1 = \lambda(0) = \lambda( px) = \lambda(x)^p. $$
This shows $\im \lambda \subseteq \nu_p$. Equality follows from the
fact $\im \lambda$ is a non-trivial subgroup of the simple abelian 
group $\nu_p$.

If $\cha F =0$ then, since $p \in \maxi$,
there exists an $n \ge 0$ such that $p^n x  \in \idl$. For such $n$,
          $$ 1 = \lambda(p^nx) = \lambda(x)^{p^n} . $$
Then $\im \lambda \subseteq
\nu_{p^\infty}$. If the inclusion were proper then there would
exist $m \ge 0$ such that $\im \lambda = \nu_{p^m}$. In this case,
if $x \in F$ then
       $$ \lambda(x) = \lambda\left( p^m \cdot \frac{ x }{p^m}\right)
                    = \lambda\left(\frac{x}{p^m} \right)^{p^m}
                       = 1 $$
since $\lambda(x/p^m)$ is a $p^m$-th root of unity. As this
would contradict the non-triviality of $\lambda$, $\im \lambda =
\nu_{p^\infty}$. 
\end{proof}
Define $\rn(\lambda)$ to be the field obtained by adjoining to $\rn$
all the character values $\lambda(x)$, $x\in F$.
Define
           $$ {\mathcal P} \simeq 
\begin{cases}  {\bf Z}/p{\bf Z}, &\text{if $\cha F =p$;}\\
       {\bf Z}_p, &\text{ if $\cha F = 0$.}
\end{cases} $$
Note that $\mathcal P$ is the topological closure of the prime ring of $F$.

\begin{lemma}
There is a canonical topological isomorphism
       $$\Gal(\rn(\lambda)/\rn) \simeq {\mathcal P}^*. $$
\end{lemma}
\begin{proof} The preceding lemma ensures that $\im \lambda$
is invariant under the action of Galois, hence restriction yields a 
homomorphism
     $$ \Gal(\rn(\lambda)/\rn) \rightarrow
          \Aut(\im \lambda) \simeq 
\begin{cases} ({\bf Z}/p{\bf Z})^*, & \text{if $\cha F =p$;} \\
       {\bf Z}_p^*, & \text{if $\cha F = 0$.}\end{cases}$$
It is readily checked that this map is an isomorphism of topological
groups. The proof is completed by appealing to the description
of $\mathcal P$ given above.
\end{proof}

The pairing
        $$ (s, t) \rightarrow \lambda(st), \qquad s,t \in F,$$
is non-degenerate and leads to an identification of $F^+$ with
its Pontryagin dual [W, II.5]. The image of $s \in F$
in the dual shall be denoted $\lambda[s]$ :
      $$\lambda[s](t) = \lambda(st), \qquad t \in F. $$
Let $\mu=dt$ be a Haar measure on $F^+$. If $\phi$ is a locally constant, complex
valued function on $F$ of compact support, the Fourier 
transform ${\mathcal F}_\lambda \phi$ is the complex valued function on $F$
defined by
    $$ {\mathcal F}_\lambda \phi(s) = \int_{F} \lambda[s](t) \phi(t) \, dt ,
      \qquad s \in F. $$
It can be shown that ${\mathcal F}_\lambda \phi$ is locally constant and has
compact support. Furthermore, the general theory of Fourier transforms
asserts the existence of a positive constant $c$, depending only
on the Haar measure $dt$, such that
       $$\left({\mathcal F}_\lambda {\mathcal F}_\lambda \phi \right) (t)
               = c \phi(-t), \qquad t \in F. $$

There is a unique Haar measure on $F^+$ for which $c=1$; it shall be 
denoted $d_\lambda t$ and will be referred to as the {\it self-dual
Haar measure associated with $\lambda$}. [W, VII.2]
\begin{lemma} If $\lambda$ has level $l$ then the associated
self-dual Haar measure is characterized by the condition
\begin{equation}\int_\O \, d_\lambda t = q^{l/2}. \end{equation}
\end{lemma} 
\begin{proof}This follows from [W, Corollary $3$, VII.2]. \end{proof}
%
\begin{cor} If $s \in F^*$ then
           $$ d_{\lambda[s]} t = \vert s \vert^{1/2} d_\lambda t. $$ 
\end{cor}
\begin{proof} Since $ \id_\lambda =  s\id_{\lambda[s]}$,  
the levels $l_1$ of $\lambda$ and $l_2$ of $\lambda[s]$
satisfy the relation $l_1 = v(s) + l_2$. Therefore, Lemma $3$ yields
      $$ \int_\O \, d_{\lambda[s]}t = 
            q^{l_2 / 2} = q^{-v(s)/2} q^{l_1/2}
                = \vert s \vert^{1/2} \int_\O \, d_\lambda t. $$
This completes the proof of the corollary. 
\end{proof}

\section{  The Schr\"odinger and Weil Representations}

Let $\aform~~$ be a non-degenerate, alternating, $F$-bilinear form
on a finite dimensional $F$-vector space $V$. The {\it Heisenberg group} 
$H$ is the group on $V \times F$ having 
multiplication
          $$(v,t)(v',t') = \left(v + v', t+t'+{\aform{v}{v'}/2}\right), 
              \qquad t,\,t' \in F,\,  v,\,v' \in V . $$

Let $\lambda$ be a non-trivial, continuous, complex linear character of $F^+$.
Since $Z(H) = 0 \times F \simeq F^+$, it may be viewed as a character
of the center of the Heisenberg group $H$.   
\begin{S-vN}There exists a smooth, irreducible
representation of $H$ having central character $\lambda$. Such a representation
is necessarily admissible, and is unique up to isomorphism.
\end{S-vN}
A proof of the Stone-von Neumann Theorem can be found in~[MVW, 
$2$.I]. The  representation provided by
the Stone-von Neumann Theorem is referred to as {\it the Schr\"odinger 
representation of type $\lambda$}.  

The symplectic group
   $$ \Sp(V) = \{\, g \in \GL(V) \, : \, \aform{vg}{wg} 
               = \aform{v}{w}, v, w \in V\, \}$$
acts on the Heisenberg group $H$ as a group of automorphisms as follows: if
$g \in \Sp(V)$ and $(t,v) \in H$ then
      $$ (t,v)g = (t, vg). $$
Given a Schr\"odinger representation $S_\lambda$ of type $\lambda$
and $g \in \Sp(V)$,
consider the representation $S_\lambda^g$ of $H$ defined by
         $$ S_\lambda^g(h) = S_\lambda(hg), \qquad 
                h \in H. $$
It is readily verified that $S_\lambda^g$ is a smooth, irreducible
representation of $H$. Furthermore, observing that $g$ acts trivially 
on $Z(H)$, $S_\lambda^g$ has central character $\lambda$.  The Stone-von
Neumann Theorem allows us to conclude that the representation $S_\lambda$
and $S_\lambda^g$ are equivalent, hence the ambient space affording $S_\lambda$ 
admits an operator $W_\lambda(g)$ for which 
     $$ S_\lambda^g(h) = W_\lambda(g)^{-1} S_\lambda(h) W_\lambda(g), 
        \qquad h \in H$$
In light of Schur's Lemma, the operator $W_\lambda(g)$ is uniquely defined 
up to multiplication by a non-zero constant. As a result, the map
        $$ g \mapsto W_\lambda(g), \qquad g \in \Sp(V), $$
is a projective representation of $\Sp(V)$, called a {\it Weil representation
of type $\lambda$}.

In this paper we will consider the Schr\"odinger models of $S_\lambda$
and $W_\lambda$ ([K, Lemma 2.2, Proposition 2.3], 
[MVW, 2.I.4(a), 2.II.6], [RR, \S $3$]). Let 
         $$V = X + Y$$
where $X$ and $Y$ are maximal, totally isotropic subspaces. The 
Schr\"odinger model is realized in the Bruhat-Schwartz 
space ${\mathcal S}(X)$ of locally constant 
functions $f : X \rightarrow \cn$  of compact support:  
if $x \in X$, $y \in Y$ and $t \in F$
then $S_\lambda((x+y,t))$ is the operator defined by 
     $$\left[S_\lambda((x+y,t))\phi\right](x') 
   = \lambda\left( t + \frac{\aform{x}{y} }{2}
                               + \aform{x'}{y}\right) \phi(x+x'), \qquad 
        \phi \in {\mathcal S}(X), x' \in X. $$

The description of the Weil representation requires some additional
notation. Viewing $x+y \in V$ as a row vector $(x,y)$, each $g \in \Sp(V)$
can be expressed in the matrix form 
\begin{equation}
       g= \begin{pmatrix} a & b\\
                          c & d \end{pmatrix} , \end{equation}
where $a : X \rightarrow X$, $b : X \rightarrow Y$, $c : Y \rightarrow X$,
and $d : Y \rightarrow Y$. With this notation,  set
         $$ Y_g =  Y / \ker c . $$
If $\mu_g$ is a Haar measure on $Y_g$ then the action of $W_\lambda(g)$
on ${\mathcal S}(X)$ is given by
\begin{multline}
       \left[ W_\lambda(g) \phi \right] (x) 
               = \\
\int_{Y_g} \left[ \lambda\left(\frac{
    \aform{xa}{xb} - 2\aform{xb}{yc} + \aform{yc}{yd}}{2}\right) 
                       \phi(xa + yc) \right] \, d\mu_g y, \quad
             \phi \in {\mathcal S}(X), \quad x \in X. 
\end{multline}
Note that the integral appearing in (3) is well-defined,
for the integrand is constant on the cosets of $\ker c$,
hence can be viewed as a function on $Y_g$. The fact $\phi \in {\mathcal S}(X)$
can be used to show that the integrand belongs to ${\mathcal S}(Y_g)$, hence
the integral converges, and that the resulting function $W_\lambda(g) \phi$
belongs to ${\mathcal S}(X)$.
\par

We now recall a particular choice of Haar measures $\mu_{\lambda,g}$ 
on $Y_g$, $g \in \Sp(V)$ [RR, \S 3.3]. Fix a 
basis $x_1, \dots, x_n$ of $X$ and let $y_1,\ldots, y_n$ be the
dual basis of $Y$ defined by the conditions
 $$ \aform{x_i}{y_j} = \delta_{ij}, \qquad 1 \le i,j \le n.$$ 
Let $\tau_i$, $0 \le i \le n$, be the element of $\Sp(V)$ defined
by
$$  x_j \tau_i = \begin{cases} -y_j, & \text{if $j \le i$;} \\
                           x_j, & \text{if $i < j$.} \end{cases} 
\quad \text{ and } 
      \quad
     y_j \tau_i = \begin{cases} x_j, & \text{if $j \le i$;}\\ 
                          y_j, & \text{if $i < j $.}\end{cases} $$
We note that $Y_{\tau_i}$ can be identified with the subspace
of $Y$ spanned by the elements $y_1,\ldots,y_i$. We define
\begin{equation}
       d\mu_{\lambda,\tau_i} y = \prod_{k=1}^i d_\lambda y_k. 
\end{equation}
where $d_\lambda y_k$ is the self-dual Haar measure associated
with $\lambda$.

Let
          $$ P = \{ g \in \Sp(V) \, : \, Yg = g \}, $$
the parabolic subgroup that leaves $Y$ invariant. If $\dim Y_g = i$ 
then [RR, Theorem 2.14] ensures the existence of elements $p_1$ 
and $p_2$ of $P$ such that
\begin{equation}
          g = p_1 \tau_i p_2. \end{equation}
Observing that the operator $p_1$ induces an isomorphism $\overline{p_1}
: Y_g \rightarrow Y_{\tau_i}$, we set
\begin{equation}  \mu_{\lambda,g} = \vert \det(p_1p_2 \vert_Y) \vert^{-1/2} \, 
           \overline{p_1} \cdot \mu_{\lambda,\tau_i}. \end{equation}
Here, $\overline{p_1} \cdot \mu_{\lambda, \tau_i}$ denotes
the pullback of the Haar measure $\mu_{\lambda,\tau_i}$ to $Y_g$
via $\overline{p_1}$ : if $E$ is a measurable subset of $Y_g$ then
        $$ \overline{p_1} \cdot \mu_{\lambda, \tau_i} (O) = 
            \mu_{\lambda, \tau_i} ( O\overline{p_1} ). $$
\begin{theorem} The measures $\mu_{\lambda,g}$, $ g \in Sp(V)$,
are well-defined. Furthermore, the projective representation $W_\lambda$
of $\Sp(V)$ defined by (3) with the Haar measures $\mu_g = \mu_{\lambda,g}$
has the following properties. 
\item{(i)} If $g \in \Sp(V)$ and $p_1$, $p_2 \in P$ then
$W_\lambda(p_1gp_2) = W_\lambda(p_1) W_\lambda(g) W_\lambda(p_2)$; in 
particular $W_\lambda$ restricts to an ordinary representation of $P$.
\item{(ii)} If $\phi \in {\mathcal S}(X)$ and
             $ p = \begin{pmatrix} a & b \\
                                    0 & d \end{pmatrix} \in P $
then
       $$ \left[W_\lambda(p)\phi\right](x) =
              \vert \det a \vert^{1/2} \lambda\left(\frac{\aform{xa}{xb}} {2}
              \right) \phi(xa), \qquad x \in X. $$
\end{theorem}

\begin{proof} This follows from
 [RR, Theorem 3.5] \end{proof}

\begin{lemma} If $s \in F^*$ and $g \in \Sp(V)$ then $\mu_{\lambda[s],g} 
= \vert s \vert_{Y_g}^{1/2} \mu_{\lambda,g}$.
\end{lemma}
\begin{proof} In light of the Corollary to Lemma 3, (4) yields
       $$ d\mu_{\lambda[s], \tau_i} y =
              \prod_{k=1}^i d_{\lambda[s]} y_k =
           \prod_{k=1}^i  \left[ \vert s \vert^{1/2} d_\lambda y_k \right]
         = \vert s \vert^{i/2} \left[ \prod_{k=1}^i d_\lambda y_k \right]
            = \vert s \vert^{i/2} d\mu_{\lambda,\tau_i}y. $$
Therefore, $(6)$ gives
      $$ \mu_{\lambda[s],g} = \left\vert \det(p_1p_2 \vert_Y) 
      \right\vert^{-1/2} \, 
           \overline{p_1} \cdot \mu_{\lambda[s], \tau_i} =
         \vert s \vert^{i/2} \left\vert \det(p_1p_2 \vert_Y) \right\vert^{-1/2}
	  \, 
           \overline{p_1} \cdot \mu_{\lambda, \tau_i}
            = \vert s \vert^{i/2} \mu_{\lambda, g} 
              = \vert s \vert^{1/2}_{Y_g}, \mu_{\lambda,g}, $$
since $Y_g$ has dimension $i$ over $F$.
\end{proof}

\section{ Action of Symplectic Similitudes}
Given $s \in F^*$, let $f_s$ be the element of $\GL(V)$ defined by
           $$ (x+y) f_s = x + sy, \qquad x \in X, y \in Y. $$
Conjugation by $f_s$ leaves the symplectic group $\Sp(V)$
invariant. In fact, if $g \in \Sp(V)$ is expressed in the matrix
form (2) then
\begin{equation}
       g^{f_s} = \begin{pmatrix} a & sb \\
                           s^{-1}c & d \end{pmatrix}. 
\end{equation}
In particular, we note that the spaces $Y_g$ and $Y_{g^{f_s}}$ are equal,
since $\ker c = \ker s^{-1}c$.
\begin{lemma}If $s \in F^*$ then $\mu_{\lambda, g^{f_s}}
= \vert s \vert^{-1/2}_{Y_g} \mu_{\lambda,g}$.
\end{lemma}
\begin{proof} Let $p_{i,s}$, $0 \le i \le n$, be the elements
of $\Sp(V)$ defined by
 $$  x_j p_{i,s} = \begin{cases} s^{-1}x_j, & \text{if $j \le i$;} \\
                           x_j, & \text{if $i < j$.} \end{cases}
 \quad \hbox{ and }   \quad
     y_j p_{i,s}= \begin{cases} sy_j, & \text{if $j \le i$;} \\
                          y_j, & \text{if $i < j $.} \end{cases}  $$
Note that $p_{i,s} \in P$ and
          $$ \det (p_{i,s}\vert_Y)  = s^i. $$
Moreover, one readily verifies that
         $$ \tau_i^{f_s} = \tau_i p_{i,s}. $$

Let $g \in G$. If $g = p_1 \tau_i p_2$, $p_1$, $p_2 \in P$, then
       $$ g^{f_s} = \left( p_1 \tau_i p_2\right)^{f_s}
                  = p_1^{f_s} \tau_i^{f_s} p_2^{f_s}
                  = p_1^{f_s} \tau_i (p_{i,s}p_2^{f_s}). $$
Observing that both $p_1^{f_s}$ and $p_{i,s}p_2^{f_s}$ belong to $P$,
$(6)$ yields
      $$ \mu_{\lambda,g^{f_s}} 
        = \vert \det (p_1^{f_s} p_{i,s}p_2^{f_s}\vert_Y) 
            \vert^{-1/2} \overline{p_1^{f_s}} \cdot \mu_{\lambda,\tau_i}. $$
Using (7), if $p \in P$ then $p^{f_s} \vert_Y = p \vert_Y$. As a 
consequence,
    $$\overline{p_1^{f_s}} = \overline{p_1} : Y_g \rightarrow Y_{\tau_i}. $$ 
In light of these observations,
        $$ \det(p_1^{f_s}p_{i,s}p_2^{f_s} \vert_Y)
             = \det (p_1 p_{i,s}p_2\vert_Y) 
              = \det(p_{i,s}\vert_Y) \cdot \det(p_1 p_2\vert_Y)  
               = s^i \det (p_1 p_2\vert_Y), $$
hence
     $$\mu_{\lambda,g^{f_s}} = \vert s^i\det (p_1p_2\vert_Y) 
            \vert^{-1/2} \overline{p_1} \cdot \mu_{\lambda,\tau_i} = 
             \vert s \vert^{-i/2} \mu_{\lambda,g} 
           = \vert s \vert^{-1/2}_{Y_g} \mu_{\lambda,g}, $$
since $Y_g$ has dimension $i$ over $F$.
\end{proof}

Let $W_\lambda^{f_s}$ be the projective representation of $\Sp(V)$ defined
by
      $$ W_\lambda^{f_s}(g) = W_\lambda\left( g^{f_s} \right). $$
For the proof of the next result, let $\vert \alpha \vert_V$ denote
the module of an automorphism $\alpha$ of an $F$-vector space $V$
[W, I.2]. We have
$$ \vert \alpha \vert_V = \left\vert \det\alpha \right\vert. $$
In particular, the module of left multiplication by 
$s\in F^*$ on $V$ satisfies
         $$ \vert s \vert_V = \vert s \vert^{\dim V}. $$ 
\begin{prop} If $s \in F^*$ then $W_\lambda^{f_s}
     = W_{\lambda[s]}$.
\end{prop}
\begin{proof} Let $g \in \Sp(V)$. We assume that $g$
has the matrix representation (2), hence that of $g^{f_s}$
is given by (7). If $\phi \in {\mathcal S}(X)$ and $x\in X$ 
then the integral formula (3) and Lemma 6 yield
   $$\begin{aligned} \left[ W_\lambda(g^{f_s})\phi\right](x) 
        &= \int_{Y_{g^{f_s}}} \left[ \lambda\left( \frac{ \aform{xa}{sxb} 
              -2 \aform{sxb}{s^{-1}yc} + \aform{s^{-1}yc}{yd} }{2}\right)
           \phi(xa + s^{-1}yc) \right] \, d\mu_{\lambda, g^{f_s}}y \\
        &= \vert s \vert_{Y_g}^{-1/2} \int_{Y_g} \left[ \lambda\left( \frac
             { \aform{xa}{sxb} 
              -2 \aform{sxb}{s^{-1}yc} + \aform{s^{-1}yc}{yd} }{2}\right)
           \phi(xa + s^{-1}yc) \right] \, d\mu_{\lambda,g}y. 
\end{aligned}  $$
Replacing $y$ by $sy$, the definition of $\vert s \vert_{Y_g}$ and Lemma $4$
yield
   $$\begin{aligned}\left[W_\lambda(g^{f_s})\phi\right](x) 
    &=\vert s \vert_{Y_g}^{-1/2} \vert s \vert_{Y_g}
             \int_{Y_g} \left[ \lambda\left(\frac{\aform{xa}{sxb} 
               -2 \aform{sxb}{yc} + \aform{yc}{syd} }{2}\right)
              \phi(xa + yc) \right] \, d\mu_{\lambda g} y  \\
        &=\vert s \vert_{Y_g}^{1/2} \int_{Y_g} 
          \left[ \lambda\left(s \cdot\frac{\aform{xa}{xb} 
               -2 \aform{xb}{yc} + \aform{yc}{yd} }{2}\right)
              \phi(xa + yc) \right] \, d\mu_{\lambda,g} y\\
        &= \vert s \vert_{Y_g}^{1/2} \int_{Y_g} 
          \left[ \lambda[s] \left( \frac{\aform{xa}{xb} 
               -2 \aform{xb}{yc} + \aform{yc}{yd} }{2}\right)
              \phi(xa + yc) \right] \, d\mu_{\lambda,g} y\\ 
        &=\int_{Y_g} 
          \left[ \lambda[s] \left(\frac {\aform{xa}{xb} 
               -2 \aform{xb}{yc} + \aform{yc}{yd} }{2}\right)
              \phi(xa + yc) \right] \, d\mu_{\lambda[s],g}y \\
        &= [W_{\lambda[s]}(g)\phi](x). 
\end{aligned} $$
This completes the proof of the proposition. 
\end{proof}

\section{Action of Galois}

Let $\mu$ be a Haar measure on a totally disconnected topological group $A$.
If $O_1$ and $O_2$ are non-empty compact open sets in $A$ then the ratio
       $$ (O_1 :  O_2) = \frac{\mu (O_1) }{\mu (O_2)} $$
is a rational number [C, I.1.1.]. Hence, if $\mu(O)$ lies in a subfield $L$
of ${\bf C}$ for some non-empty compact open set $O$ then the same is
true for all non-empty compact open sets. The measure $\mu$ is
said to~{\it $L$-rational} if this is the case.

\begin{lemma} The measures $\mu_{\lambda, g}$, $g \in \Sp(V)$,
are ${\bf Q}(\sqrt{q})$-rational.
\end{lemma}
\begin{proof}
If $t \in F^*$ then $\vert t \vert$ is
a power of $q$. Therefore, (6) shows that it is sufficient to
verify that the measures $\mu_{\lambda, \tau_i}$ 
are ${\bf Q}(\sqrt{q})$-rational. Formulas $(1)$ and $(4)$
ensure that this is indeed the case: if ${\mathcal Y}_i 
= \sum_{k=1}^i \O y_k$ then 
        $$ \int_{{\mathcal Y}_i} \, d\mu_{\lambda,\tau_i} y = q^{il/2}. $$
This completes the proof of the lemma.  \end{proof}

Let $A$ be a totally disconnected topological group. If $L$ is a subfield
of ${\bf C}$, let ${\mathcal S}(A,L)$ denote the space of locally constant,
$L$-valued functions on $A$ of compact support. 
\begin{lemma} 
 Let $A$ be a totally disconnected topological group,
$L \subseteq K$ an extension of fields, and $\mu$ a $L$-rational
Haar measure on $A$. If $\phi \in {\mathcal S}(A,K)$ then $\int_A \phi \, d\mu$
belongs to $K$. \end{lemma}
\begin{proof} Since $\phi \in {\mathcal S}(A,K)$, there exists 
compact open subsets $A_1$,$\ldots$, $A_k$ of $A$ and 
scalars $c_1$, $\ldots$, $c_k$ in $K$ such that
       $$ \phi = \sum_{i=1}^k c_i \chi_{A_i}. $$
Here, $\chi_{A_i}$ denotes the characteristic function of $A_i$. 
Since $\mu(A_i) \in L \subset K$, it follows that
    $$\int_A \phi \, d\mu = \sum_{i=1}^k c_i \mu(A_i) $$
lies in $K$.
\end{proof}
  
Let $\rn(\lambda)$ be the character field of $\lambda$ and
set
               $$ E = {\bf Q}(\lambda)(\sqrt{-1}). $$
Observe that Lemma $1$ ensures that ${\bf Q}(\sqrt{q})$ is a subfield of $E$.
\begin{prop} The operators $W_\lambda(g)$, $g \in \Sp(V)$, leave
the subspace ${\mathcal S}(X, E)$ invariant.
\end{prop}
\begin{proof}  
If $\phi \in {\mathcal S}(X, E)$ then the integrand
in (3) lies in ${\mathcal S}(Y_g, E)$, since ${\bf Q}(\lambda) \subseteq E$. In 
light of Lemma $8$, Lemma $9$ applied in
the case $A = Y_g$, $K = E$, $L = {\bf Q}(\sqrt{q})$, 
and $\mu = \mu_{\lambda,g}$ allows us to deduce
that the integral (3) lies in $E$. It follows immediately 
that $W_\lambda(g)\phi \in {\mathcal S}(X, E)$. \end{proof}

By Lemma 1, $E$ is a Galois extension of ${\bf Q}$.  Its Galois 
group acts on ${\mathcal S}(X, E)$ : if $\sigma \in \Gal(E/{\bf Q})$ 
and $\phi \in {\mathcal S}(X,E)$ then
\begin{equation}
     (\sigma(\phi))(x) = \sigma(\phi(x)), \qquad x \in X.
\end{equation}
There is an associated Galois action on $\End{\mathcal S}(X,E)$ : if $\sigma
\in G$ and $T \in \End {\mathcal S}(X,E)$ then 
\begin{equation}
        {}^\sigma T  (\phi) 
               = \sigma\left[ T\left({\sigma^{-1}}(\phi) 
          \right) \right] , \qquad \phi \in {\mathcal S}(X,E). 
\end{equation}
The Galois group also permutes the unitary characters of $F^+$: 
if $\sigma \in \Gal(E/\rn)$ and $\lambda$
is a unitary character of $F^+$ then $^\sigma\lambda$
is the character defined by
        $$ ^\sigma\lambda(t) = \sigma(\lambda(t)), \qquad
              t \in F^+. $$

Let $\mathcal P$ be the topological closure of the prime ring of $F$. 
The image of $s \in {\mathcal P}^*$ in $\Gal(\rn(\lambda)/\lambda)$ 
under the canonical isomorphism of Lemma $2$ will be 
denoted $\sigma_s$. 
\begin{lemma} 
Let $\sigma \in \Gal(E/\rn)$. 
If $\sigma\vert_{\rn(\lambda)} = \sigma_s$ 
then $^\sigma\lambda = \lambda[s] $.
\end{lemma} 
\begin{proof} ($\cha F = 0$) Let $\id$ be the conductor 
of $\lambda$. Given $t \in F$, fix $n \ge 1$ such that $t \in p^{-n}\id$.
Since $p^n t \in \id$,
            $$ 1 = \lambda(p^n t) = \lambda(t)^{p^n} , $$
thus $\lambda(t) \in \nu_{p^n}$. Fixing $r \in {\bf Z}$ such that $s \equiv r
\bmod p^n {\mathcal P}$, 
      $$ \left(^\sigma\lambda \right) (t)
           = \sigma\left( \lambda(t) \right) 
           = \lambda(t)^r = \lambda(rt) = \lambda(st), $$
the last equality following from the fact $rt \equiv st \bmod \id$. \end{proof}

Given $\sigma \in \Gal(E/\rn)$, let $^{\sigma}W_\lambda$ be the
projective representation defined by
       $$ \left(^{\sigma}W_\lambda\right)(g) = {^\sigma(W_\lambda(g))},
          \qquad  g \in \Sp(V).$$
\begin{prop}Let $\sigma \in {\rm Gal}(E/{\bf Q}(\sqrt{q}))$.
If $\sigma \vert_{\rn(\lambda)} = \sigma_s$ 
then $^{\sigma}W_\lambda(g) = W_{\lambda[s]}(g)$.
\end{prop} 

The proof of Proposition 12 is based on the integral formula (3)
and the following 
\begin{lemma} 
Let $A$ be a totally disconnected topological group,
$L \subseteq K$ an extension of fields, and $\mu$ a $L$-rational Haar
measure on $A$. If $\sigma$ is an $L$-automorphism of $K$ then,
for all $\phi \in {\mathcal S}(A, K)$, 
       $$ \int_A \sigma(\phi) \, d\mu =
              \sigma\left( \int_A \phi \, d\mu \right). $$
\end{lemma}
\begin{proof} Using the notation introduced
in the proof of Lemma $5$, if
        $$ \phi = \sum_{i=1}^k c_i \chi_{A_i} $$
then
        $$ \sigma(\phi) = \sum_{i=1}^k \sigma(c_i) \chi_{A_i}. $$
Therefore, since $\mu(A_i) \in L$ is fixed by $\sigma$,
    $$ \int \sigma(\phi) \, d\mu =
                   \sum_{i=1}^k \sigma (c_i) \mu(A_i) 
             = \sum_{i=1}^k \sigma(c_i) \sigma (\mu(A_i)) 
             = \sigma \left( \sum_{i=1}^k c_i \mu(A_i) \right) 
                = \sigma\left( \int_{A} \phi \, d\mu \right). $$
This completes the proof of the lemma.  \end{proof} 

\begin{proof}[Proof of Proposition 12]
Let $g \in 
\Sp(V)$, $\phi \in {\mathcal S}(X,E)$, and $x \in X$. We assume
$g$ has the matrix representation (2). Lemma $8$ asserts 
that the measure $\mu_{\lambda,_g}$  is ${\bf Q}(\sqrt{q})$-rational. Applying 
Lemma 13 to the case $A = Y_g$, $L = {\bf Q}(\sqrt{q})$, $K = E$, 
and $\mu = \mu_{\lambda,g}$, the definition of ${}^\sigma \! W_\lambda$,
the formula (3), and Lemma 11 yield
  $$ \begin{aligned} \left[ ^{\sigma}W_\lambda(g) \phi \right](x) 
         &= \sigma \left[ W_\lambda(g) (\sigma^{-1}\phi)(x) 
              \right]\\
         &= \sigma \left[ \int_{Y_g} \left[ \lambda\left(\frac {\aform{xa}{xb} - 2
                    \aform{xb}{yc} + \aform{yc}{yd} }{2 }\right) 
                       (\sigma^{-1}\phi)(xa + yc) \right] \, 
                             d\mu_{\lambda,g}y \right] \\
         &=\int_{Y_g} \!\left[ {^\sigma\lambda}\left(\frac {\aform{xa}{xb} - 2
                    \aform{xb}{yc} + \aform{yc}{yd} }{2 }\right) 
                      \phi(xa + yc) \right] \, d\mu_{\lambda,g}y \cr
         &= \int_{Y_g} \!\left[ \lambda[s]\left(\frac {\aform{xa}{xb} - 2
                    \aform{xb}{yc} + \aform{yc}{yd} }{2 }\right) 
                      \phi(xa + yc) \right] \, d\mu_{\lambda,g}y. 
                   \end{aligned} $$   
Observing $s \in {\mathcal P}^* \subseteq \O^*$, Lemma 5 
implies that $\mu_{\lambda[s],g}= \mu_{\lambda,g}$. The preceding calculation 
thus gives
    $$ \begin{aligned} \left[ ^\sigma{W_\lambda}(g) \phi \right](x) 
             &= \int_{Y_g} \!\left[ \lambda[s]\left(\frac {\aform{xa}{xb} - 2
                     \aform{xb}{yc} + \aform{yc}{yd} }{2 }\right) 
                       \phi(xa + yc) \right] \, d\mu_{\lambda[s],g}y \\
            &=  \left[W_{\lambda[s]}(g) \phi \right](x). \end{aligned}$$
This completes the proof. \end{proof} 

\section{The Fundamental Identity}
Let
       $$ {\mathfrak G} = \{ \sigma \in \Gal(E/\rn(\sqrt{q})) \, : \,
             \exists s \in \O^* \hbox{ such that } 
             \sigma \vert_{\rn(\lambda)} = \sigma_{s^2} \}. $$
Note that ${\mathfrak G}$ is a subgroup of $\Gal(E/\rn(\sqrt{q}))$.
Given $s \in F^*$, let $g_s \in \Sp(V)$ be the map defined by
           $$ (x+y) g_s = s^{-1} x + sy, \qquad x \in X, y \in Y. $$
We observe that $g_s$ lies in the parabolic subgroup $P$ that leaves $Y$
invariant and is related to the operator $f_{s^2}$ defined earlier 
by the identity
        $$ f_{s^2} = sI \circ g_s.$$
\begin{prop} 
 Let $\sigma \in {\mathfrak G}$ and $g \in \Sp(V)$. 
If $\sigma \vert_{\rn(\lambda)} = \sigma_{s^2}$, $s \in \O^*$, then
       $$ ^{\sigma}W_\lambda(g) =
             W_\lambda(g_s)^{-1} W_\lambda(g) W_\lambda(g_s). $$
\end{prop}
\begin{proof} In light of Propositions $7$ and $12$,
    $$ ^{\sigma}W_\lambda(g) = W_{\lambda[s^2]}(g)
                = W_\lambda^{f_{{s^2}}} (g)
                   = W_\lambda( g^{f_{s^2}} )
          = W_\lambda( g^{g_s} ). $$
Applying Theorem 4(i) with $p_1^{-1} = p_2 = g_s$,
          $$ W_\lambda(g^{g_s}) 
             = W_\lambda(g_s^{-1})W_\lambda(g)W_\lambda(g_s)
             = W_\lambda(g_s)^{-1} W_\lambda(g) W_\lambda(g_s). $$
This completes the proof of the proposition. \end{proof}
\begin{cor}If $t \in F^*$ and $\sigma \in {\mathfrak G}$ 
 then $^{\sigma}W_\lambda(g_t) = W_\lambda(g_t).$ 
\end{cor}
\begin{proof} Fix $s \in \O^*$ such 
that $\sigma \vert_{\rn(\lambda)} = \sigma_{s^2}$. Observing 
that $g_s$ and $g_t$ are commuting elements of $P$, the preceding proposition 
combines with Theorem 4(i) to yield
  $$ ^{\sigma}W_\lambda(g_t) = W_\lambda(g_s)^{-1} W_\lambda(g_t)
                        W_\lambda(g_s) 
      = W_\lambda(g_s^{-1}g_tg_s)  = W_\lambda(g_t), $$ 
as required. \end{proof}
   
\section{The Cocycle}
Define
     $$ {\mathfrak H} = \Gal\left(E/\rn(\sqrt{p},\sqrt{-p})\right). $$
Our aim in this section is the construction of a $1$-cocycle
$\delta : {\mathfrak H} \rightarrow \GL\left({\mathcal S}(X,E)\right)$
such that
\begin{equation}
       {}^{\sigma}W_\lambda(g) = \delta(\sigma)^{-1}W_\lambda(g)
                           \delta(\sigma), \qquad g \in \Sp(V),\, 
             \sigma \in {\mathfrak H}. 
\end{equation}
Let $p^* = (-1)^{(p-1)/2}p$. When combined with restriction to $\rn(\lambda)$,
the canonical isomorphism of Lemma $2$ yields 
\begin{equation}
        {\mathfrak H} \simeq \Gal\left(\rn(\lambda)/ \rn(\sqrt{p^*})\right)
       \simeq \left({\mathcal P}^*\right)^2. 
\end{equation}
If $\epsilon \in {\mathcal P}^*$ is a primitive $p-1$ root of unity
then
       $$ {\mathcal P}^* = \left< \epsilon \right> \times U_1 $$
with
          $$ U_1 = \begin{cases} \{ 1 \}, & \text{if $\cha F = p$;} \\
                 \{ r \in {\mathcal P} \, : \, r \equiv 1 \bmod p\}, &
                            \text{if $\cha F = 0$,}\end{cases}$$
a pro-$p$ group. As $p$ is assumed to be odd, the map $r \mapsto r^2$
is an automorphism of $U_1$, hence
        $$ {{\mathcal P}^*}^2 =  \left< \epsilon^2 \right> \times U_1. $$
The isomorphism $(11)$ identifies $U_1$ 
with $\Gal\left(E/ \rn(\nu_p, \sqrt{-1})\right)$, where $\nu_p$
is the group of complex $p$-th roots of unity. This in turn leads to 
an identification of $\left<\epsilon^2\right>$ with
 $$ {\mathfrak H} / \Gal\left(E/ 
           \rn(\nu_p, \sqrt{-1}) \right)
               \simeq \Gal\left(\rn(\nu_p, \sqrt{-1})/
             \rn(\sqrt{p}, \sqrt{-p})\right). $$
In particular, the element $\eta$ of ${\mathfrak H}$ characterized by
\begin{equation}\label{eta}
              \eta\vert_{\rn(\lambda)} = \sigma_{\epsilon^2}
\end{equation}
has order $(p-1)/2$ and restricts to a generator 
of $\Gal\left(\rn(\nu_p, \sqrt{-1})/\rn(\sqrt{p}, \sqrt{-p})\right)$.

Given $\sigma \in {\mathfrak H}$, there is a unique integer $i$, $1 \le i \le 
(p-1)/2$, and a unique element $s \in U_1$, such that
           $$ \sigma \vert_{\rn(\lambda)} = \sigma_{\epsilon^{2i}s^2}. $$
If $\tau$ is a second element of ${\mathfrak H}$, say
          $$ \tau \vert_{\rn(\lambda)} = \sigma_{\epsilon^{2j}t^2},
              \qquad 1 \le j \le (p-1)/2, \quad t \in U_1. $$           
then
        $$ \sigma\tau \vert_{\rn(\lambda)} = \sigma_{\epsilon^{2k}}(st)^2, $$
where $st \in U_1$ and
            $$ k =\begin{cases} i + j, & \text{if $i+j \le (p-1)/2$;} \\
                i + j - \frac{p-1 }{2}, & \text{if $i + j > (p-1)/2$.}
                   \end{cases}$$ 

Our initial attempt at the construction of the cocycle is to define
           $$ D(\sigma) = W_\lambda(g_{\epsilon^is}), \qquad 
              \sigma\vert_{\rn(\lambda)} = \sigma_{\epsilon^{2i}s^2}, 
                        \quad 1 \le i \le (p-1)/2, \quad s \in U_1. $$
Proposition $14$ ensures that 
\begin{equation}
         {}^{\sigma}W_\lambda(g) = D(\sigma)^{-1}W_\lambda(g)
                           D(\sigma), \qquad g \in \Sp(V), 
             \sigma \in {\mathfrak H}.
\end{equation}
Assuming $\sigma$ and $\tau$ are as above, the definition of $D$ yields
          $$ D(\sigma\tau) = W_\lambda(g_{\epsilon^kst}). $$
On the other hand, the Corollary to Proposition $14$ gives
          $$ {}^\sigma D(\tau) = {}^\sigma W_\lambda(g_{\epsilon^jt}) 
                   = W_\lambda(g_{\epsilon^jt}), $$
hence Theorem 4(i) yields
       $$ D(\sigma) {}^\sigma D(\tau)
               = W_\lambda(g_{\epsilon^is}) W_\lambda(g_{\epsilon^jt})
                   = W_\lambda(g_{\epsilon^{i+j} st}). $$
If $i+j \le (p-1)/2$ then
            $$ W_\lambda(g_{\epsilon^{i+j} st}) = W_\lambda(g_{\epsilon^k}).$$
If $i+j > (p-1)/2$ then, since $\epsilon^{(p-1)/2}= -1$, Theorem 4(i)
yields
      $$ W_\lambda(g_{\epsilon^{i+j} st}) = W_\lambda(g_{-\epsilon^k st})
            = W_\lambda(\iota) W_\lambda(g_{\epsilon^k st}), $$
where $\iota=g_{-1}$ is the central involution of $Sp(V)$ that maps $v \in V$
to $-v$. In summary,
\begin{equation}\label{old15} 
        D(\sigma) {}^\sigma\! D(\tau)  = 
            \begin{cases} D(\sigma \tau), & \text{if $i+j \le
                        (p-1)/2$;} \\
               W_\lambda(\iota) D(\sigma\tau),
                     & \text{if $i+j > (p-1)/2$.} \end{cases}
\end{equation}
In particular, $D$ is not a $1$-cocyle; to get one we must account for the
factor $W_\lambda(\iota)$.

Since $\iota \in P$, Theorem 4(i) implies that if $\phi$ belongs 
to ${\mathcal S}(X,E)$ then
         $$\left[W_\lambda(\iota)\phi\right](x) = \phi(-x), \qquad
                x\in X. $$
In particular, $W_\lambda(\iota)$ is an involution, hence
the operators
       $$ \rho_e = \frac{1 }{2} \left( I+ W_\lambda(\iota) \right) \quad
              \text{ and } \quad \rho_o =\frac {1 }{2} \left( 
                  I - W_\lambda(\iota) \right) $$
are orthogonal idempotents. Furthermore, recalling $\iota = g_{-1}$, the 
Corollary to Proposition 14 shows that both $\rho_e$ and $\rho_o$ are fixed by
the action of Galois. Finally, since $I = \rho_e + \rho_o$, it is easily
verified that the operators 
            $$ \rho_e + c \rho_o, \qquad c \in E, c \ne 0, $$
are invertible.

\begin{lemma} 
 The norm equation
       $$ N(u) = -1, \qquad N : \rn(\nu_p, \sqrt{-1}) \rightarrow 
                            \rn(\sqrt{p}, \sqrt{-p})$$
has a solution.
\end{lemma}
\begin{proof} The case $p \equiv 1 \bmod 4$ is covered by~[CMS,
Lemma 24(2)]. In the case $p \equiv 3 \bmod 4$ one has
       $$N(-1) = (-1)^{[ \rn(\nu_p, \sqrt{-1}) : \rn(\sqrt{p}, \sqrt{-p})]} =
             (-1)^{ (p-1)/ 2} = -1, $$
since $(p-1)/2$ is odd. \end{proof}

Let $u$ be a solution of the norm equation of the preceding lemma. 
Given $\sigma 
\in {\mathfrak H}$, set
      $$ A(\sigma) = \rho_e + \left(\prod_{l=0}^{k-1} \eta^l(u)\right) 
                      \rho_0 , \quad
       \text{where }      \sigma\vert_{\rn(\lambda)} = \sigma_{\epsilon^{2i}s^2},
         \quad 1 \le i \le (p-1)/2, \quad s \in U_1 $$
where $\eta$ satisfies (\ref{eta}).
The remarks preceding Lemma 15  ensure that $A(\sigma) \in \GL({\mathcal S}(X,E))$.
With the notation introduced earlier, if $\sigma$ and $\tau$ belong
to $\mathfrak H$ then
      $$ A(\sigma\tau) = \rho_e + \left(\prod_{l=0}^{k-1} \eta^l(u)\right)
                            \rho_0. $$
On the other hand, observing 
              $$ \sigma \eta^{-i} \vert_{\rn(\lambda)} 
              = \sigma_{\epsilon^{2i}s^2} \sigma_{\epsilon^2}^{-i}
                   = \sigma_{\epsilon^{2i}s^2} \sigma_{\epsilon^{-2i}}  
                = \sigma_{s^2},$$
the fact $(11)$ identifies $U_1$ 
with $\Gal\left(E / \rn(\nu_p, \sqrt{-1})\right)$ allows us to deduce
that the restrictions of $\sigma$ and $\eta^i$ to $\rn(\nu_p, \sqrt{-1})$
coincide. Therefore,
    $$\begin{aligned}{}^\sigma\!A(\tau) &= {}^\sigma \!\left[ \rho_e 
       + \left(\prod_{l=0}^{j-1} \eta^l(u)\right) \rho_0 \right] 
      = \rho_e + {}^\sigma\!\left(\prod_{l=0}^{j-1} \eta^l(u)\right)
                            \rho_0  \\
      &= \rho_e + {}^{\eta^i} \!\left(\prod_{l=0}^{j-1} \eta^l(u)\right)
                            \rho_0 
      = \rho_e + \left(\prod_{l=i}^{i+j-1} \eta^l(u)\right)
                            \rho_0, \end{aligned} $$
hence
    $$\begin{aligned} A(\sigma) {}^\sigma\!A(\tau) &= 
    \left[ \rho_e + \left(\prod_{l=0}^{i-1} \eta^l(u)\right)
                            \rho_0 \right]  
     \left[ \rho_e + \left(\prod_{l=i}^{i+j-1} \eta^l(u)\right)
                            \rho_0 \right] \\
      &= \left[ \rho_e + \left(\prod_{l=0}^{i+j-1} \eta^l(u)\right)
                            \rho_0 \right].\end{aligned} $$
If $i+j \le (p-1)/2$ then
          $$ \prod_{l=0}^{i+j-1} \eta^l(u) = \prod_{l=0}^{k-1} \eta^l(u), $$
hence
    $$ A(\sigma) {}^\sigma\!A(\tau) = A(\sigma\tau).$$
If $i+j > (p-1)/2$ then the choice of $\eta$ and $u$ yield
        $$ \prod_{l=0}^{i+j-1} \eta^l(u) =
            \left(\prod_{l=0}^{(p-3)/2} \eta^l(u)\right) 
            \left(\prod_{l={(p-1)/2}}^{i+j-1} \eta^l(u)\right)
                    =  N(u)   \prod_{l=0}^{k-1} \eta^l(u) 
                    = - \prod_{l=0}^{k-1} \eta^l(u). $$
Observing that $\rho_e = \rho_e W_\lambda(\iota)$ and $-\rho_o =
\rho_o W_\lambda(\iota)$, 
    $$ A(\sigma) {}^\sigma\!A(\tau) = \rho_e - \left(\prod_{l=0}^{k-1}
\eta^l(u)\right) \rho_0 = \left[ \rho_e + \left(\prod_{l=0}^{k-1}
  \eta^l(u)\right) \rho_0\right] W_\lambda(\iota) =
A(\sigma\tau)W_\lambda(\iota).  $$ In summary,
\begin{equation}\label{old14}
    A(\sigma){}^\sigma\!A(\tau) = 
    \begin{cases} A(\sigma\tau), & \text{if $i+j \le (p-1)/2$;} \\
            A(\sigma\tau)W_\lambda(\iota), &\text{if $i+j > (p-1)/2$.}
 \end{cases}
\end{equation}

Consider the map $\delta : {\mathfrak H} \rightarrow \GL({\mathcal S}(X,E))$ given by
         $$\delta(\sigma) =A(\sigma)D(\sigma). $$
If $\sigma, \tau \in \mathfrak H$ are as above 
     $$ \delta(\sigma) {}^\sigma\!\delta(\tau) =
           (A(\sigma)D(\sigma)) {}^\sigma\!(A(\tau) D(\tau)) 
           =  A(\sigma) D(\sigma) {}^\sigma\! A(\tau) {}^\sigma\!D(\tau). $$
By Theorem  4(1), ${}^\sigma\!A(\tau) \in E[W_\lambda(i)]$ commutes 
with $D(\sigma) = W_\lambda(g_{\epsilon^is})$, hence
     $$ A(\sigma)D(\sigma){}^\sigma\! A(\tau) {}^\sigma\!D(\tau)
           = A(\sigma) {}^\sigma\!A(\tau) D(\sigma) {}^\sigma\!D(\tau). $$
If $i+j > (p-1)/2$ then (14) and (15) yield
$$  A(\sigma) {}^\sigma\!A(\tau) D(\sigma) {}^\sigma\!D(\tau) =
     A(\sigma\tau) W_\lambda(\iota) W_\lambda(\iota) D(\sigma\tau)
         = A(\sigma\tau) D(\sigma\tau). $$
Since this is trivially true if $i+j \le (p-1)/2$, we conclude
     $$ \delta(\sigma) {}^\sigma\!\delta(\tau) = 
            A(\sigma\tau)D(\sigma\tau) = \delta(\sigma\tau). $$
This shows that $\delta$ is a $1$-cocycle. Furthermore, if $g \in \Sp(V)$ 
then Theorem 4(i) shows that $A(\sigma) \in E[W_\lambda(\iota)]$ 
commutes with $W_\lambda(g)$, hence $(12)$ yields
     $$\begin{aligned} \delta(\sigma)^{-1} W_\lambda(g) \delta(\sigma) &=
             (A(\sigma)D(\sigma))^{-1} W_\lambda(g) A(\sigma) D(\sigma) \\
             &= D(\sigma)^{-1} A(\sigma)^{-1} W_\lambda(g) A(\sigma)D(\sigma)\\
               &= D(\sigma)^{-1} W_\lambda(g) D(\sigma) \\ &=
                  {}^\sigma\! W_\lambda(g) \end{aligned} $$
which verifies that $(10)$ is satisfied.

\section{The Triviality of the Cocycle}
Let $\delta : {\mathfrak H} \rightarrow \GL({\mathcal S}(X,E))$ be the $1$-cocycle
satisfying (10) constructed above. 
\begin{lemma} 
If $\phi \in {\mathcal S}(X,E)$ then there exists
an open subgroup $\mathfrak K$ of $\mathfrak H$ such that
      $$ \delta(\sigma) \phi = \phi, \qquad \sigma \in {\mathfrak K}.$$
\end{lemma}
\begin{proof}  If $\cha F = p$ then ${\mathfrak H}$ is 
a finite discrete group, so one may take ${\mathfrak K}$ to be the
trivial subgroup. 

Assume $\cha F = 0$. If $\mathfrak X$ is a lattice in $X$ then the
subgroups
           $$ p^k {\mathfrak X}, \qquad k \in {\bf Z}, $$
form a local base at the origin. Therefore, given $x \in X$,
there exist $i_x \in {\bf Z}$ such that $\phi$ is constant on the 
coset $x + p^{i_x} \mathfrak X$. As the family
         $\left\{ x + p^{i_x} \mathfrak X \, : \, x \in X \right\}$
is an open cover of $X$, there exists $x_1$, $\ldots$, $x_m$
in $X$ such that
          $$\supp \phi \subseteq \bigcup_{j=1}^m {x_j + p^{i_{x_j}}{\mathfrak X}}. $$
Set
          $$ i = \max \left\{ i_{x_1}, \ldots, i_{x_m} \right\} $$
and consider ${x + p^i {\mathfrak X}} \cap \supp \phi$, $x \in X$.
If it is empty then the restriction of $\phi$ to
the coset $x+ p^i {\mathfrak X}$ is identically $0$. If not, 
there exists $j$ such that ${x + p^i {\mathfrak X}} \cap  {x_j + p^{i_{x_j}}
{\mathfrak X}}$ is non-empty, hence
     $${x + p^i {\mathfrak X}} \subseteq  {x_j + p^{i_{x_j}} {\mathfrak X}} $$
by choice of $i$. The choice of $i_{x_j}$ thus ensures that 
the restriction of $\phi$ to $x + p^i{\mathfrak X}$ is the constant function
with value $\phi(x_j)$. We conclude that $\phi$ is constant on 
the $p^i{\mathfrak X}$-cosets of $X$.

Let $\sigma \in {\mathfrak H}$. If $\sigma \vert_{\rn(\lambda)} = \sigma_{r^2}$,
$r \in U_1$, then by construction $ \delta(\sigma) = W_\lambda(g_r)$.
Observing
            $$ g_r = \begin{pmatrix} r^{-1} \cdot 1_X & 0 \\
                                  0  &  r \cdot 1_Y \end{pmatrix} 
\in P,$$
if $x \in X$ then Theorem 4(i) yields
  $$ (\delta(\sigma)\phi)(x) = 
         (W_\lambda(g_r)\phi)(x) =
            \vert r \vert^{-\dim X /2} \lambda\left( \frac{\aform{r^{-1}x}{rx}}{2} 
                  \right) \phi(r^{-1}x) 
       = \phi(r^{-1}x), $$
since $r$ is a unit and $\aform{~}{~}$ is $F$-bilinear and alternating.
Fix $j \in {\bf Z}$ such that $i > j$ and
          $$ \supp \phi \subseteq p^j {\mathfrak X}. $$
If $x \not\in p^j{\mathcal X}$ then neither is $r^{-1}x$, so the choice 
of $j$ ensures that
    $$ (\delta(\sigma)\phi)(x) = \phi(r^{-1}x) = 0 = \phi(x). $$
On the other hand, suppose $x \in p^j{\mathcal X}$. In this case, 
if $r \equiv 1 \bmod p^{i-j}$ then
           $$ r^{-1}x + p^i {\mathcal X} = x + p^{i-j}p^j{\mathfrak X} + p^i{\mathcal X}
                          = x + p^i{\mathcal X},$$
hence the choice of $i$ ensures that
        $$(\delta(\sigma)\phi)(x) = \phi(r^{-1}x) =\phi(x). $$
In light of the preceding discussion, 
        $$ {\mathfrak K} = \left\{ \sigma \in {\mathfrak H} \, : \,
               \sigma \vert_{\rn(\lambda)} = \sigma_{r^2}, 
                 r \equiv 1 \bmod p^{i-j} \right\}
          = \Gal\left(E/\rn\left(\nu_{p^{i-j}}, \sqrt{-1}\right)\right) $$
has the required properties. 
\end{proof}

Let $K/k$ be a Galois extension and $M$ a $K$-vector space equipped with 
an semi-linear action of the Galois group $\Gal(K/k)$ : 
if $\sigma \in \Gal(K/k)$, $m \in M$  and $e \in K$ then
       $$ \sigma (em) = \sigma(e) \sigma(m). $$
For such an action, the fixed-point set
   $$ M^{\Gal(K/k)} = \left\{ m \in M \, : \, m = \sigma(m) \, \forall
         \sigma \in \Gal(K/k) \right\} $$
is a $k$-vector space. The canonical action
of $\Gal(K/k)$ on $K$ yields a semi-linear action on the tensor product $K
\otimes _k M^{\Gal(K/k)}$ :
           $$ \sigma(e \otimes m) = \sigma(e) \otimes m, \qquad
                   \sigma \in \Gal(K/k), e \in E, m \in M^{\Gal(K/k)}. $$ 
The action of Galois on $M$ is said to be {\it smooth} if the 
stabilizer of each $m \in M$ is open in $\Gal(K/k)$.  
\begin{prop}[Galois Descent] If $M$ is a $K$-vector space
equipped with a semi-linear, smooth action of $\Gal(K/k)$ then
the canonical map
        $$ \psi : K \otimes_k M_k  \rightarrow M $$
is a $K$-linear isomorphism of $\Gal(K/k)$-modules.
\end{prop}
\begin{proof}
The case $K = k_s$, the separable closure
of $k$, is proved in [B, AG.14.2]. The general case is proved
using the same argument, {\it mutatis mutandis}.
\end{proof}

\begin{prop}
There exists $\alpha 
\in \GL\left({\mathcal S}(X,E)\right)$ such that
\begin{equation}
       \delta(\sigma) = \alpha^{-1} {^\sigma\alpha}, \qquad
             \sigma \in {\mathfrak H}.
\end{equation}
\end{prop}
\begin{proof}  The canonical action $(8)$ of $\mathfrak H$
on ${\mathcal S}(X,E)$ is clearly semi-linear. It is furthermore
smooth, since each element of ${\mathcal S}(X,E)$ takes only
finitely many values in $E$.

On the other hand, since $\delta$ is a 1-cocyle, then
      $$ (\sigma , \phi) \mapsto \delta(\sigma) \sigma(\phi), 
            \qquad \sigma \in {\mathfrak H}, \phi \in {\mathcal S}(X,E), $$
is also an action of $\mathfrak H$ on ${\mathcal S}(X,E)$, referred to
as the twisted action by $\delta$. It is semi-linear,
since $\delta$ takes values in $\GL\left({\mathcal S}(X,E)\right)$.
Since the original action is smooth,
if $\phi \in {\mathcal S}(X,E)$ then there exists an open subgroup ${\mathfrak H}_1$
such that
          $$ \sigma(\phi) = \phi, \qquad \sigma \in {\mathfrak H}_1. $$
Furthermore, Lemma $16$ asserts that there is an open subgroup ${\mathfrak K}$
of $\mathfrak H$ such that
          $$ \delta(\sigma)\phi = \phi, \qquad \sigma \in {\mathfrak K}. $$
Therefore, if $\sigma \in {\mathfrak H}_1 \cap {\mathfrak K}$ then
         $$ \delta(\sigma) \sigma(\phi) = \delta(\sigma) \phi = \phi. $$
This shows that the stabilizer of $\phi$ under the twisted action
contains the open subgroup ${\mathfrak H}_1 \cap {\mathfrak K}$. Since
it is the union of its ${\mathfrak H}_1 \cap {\mathfrak K}$-cosets, it
follows that the stabilzer of $\phi$ under the twisted action is open.
We conclude that the twisted action is smooth.

Using ${\mathcal S}(X,E)$ and $_\delta{\mathcal S}(X,E)$ to denote the $\mathfrak H$-modules
defined by the natural and twisted actions, respectively, Galois Descent
asserts the existence of $E$-linear, $\mathfrak H$-equivariant isomorphisms
 $$ _\delta{\mathcal S}(X,E)  \simeq  E \otimes_{\rn(\sqrt{p}, \sqrt{-p})} {
    _\delta{\mathcal S}(X,E)}^{\mathfrak H}   
\quad  \hbox{ and } \quad 
     E \otimes_{\rn(\sqrt{p}, \sqrt{-p})} {{\mathcal S}(X,E)}^{\mathfrak H}  
         \simeq  {\mathcal S}(X,E). $$ 
In particular,
    $$ \dim_{\rn(\sqrt{p},\sqrt{-p})} {_\delta{\mathcal S}(X,E)^{\mathfrak H}}    
        = \dim_E {\mathcal S}(X,E) = 
      \dim_{\rn(\sqrt{p},\sqrt{-p})} {\mathcal S}(X,E)^{\mathfrak H}, $$
so $_\delta{\mathcal S}(X,E)^{\mathfrak H}$ and ${\mathcal S}(X,E)^{\mathfrak H}$
are $\rn(\sqrt{p}, \sqrt{-p})$-isomorphic. As any such isomorphism
extends by scalars to a E-linear, $\mathfrak H$-equivariant isomorphism
       $$ E \otimes_{\rn(\sqrt{p}, \sqrt{-p})}
    {_\delta{\mathcal S}(X,E)^{\mathfrak H}}   \simeq 
          E \otimes_{\rn(\sqrt{p}, \sqrt{-p})} {{\mathcal S}(X,E)^{\mathfrak H}}  , $$
we conclude
     $$ _\delta{\mathcal S}(X,E) \simeq  {\mathcal S}(X,E). $$

Let $\alpha \in \GL\left({\mathcal S}(X,E)\right)$ be a $\mathfrak H$-equivariant
isomorphism ${_\delta{\mathcal S}(X,E)} \rightarrow {\mathcal S}(X,E)$.
If $\sigma \in {\mathfrak H}$ and $\phi \in {\mathfrak H}$ then the definition
of the twisted action ensures that
    $$ \alpha \delta(\sigma) \sigma(\phi) = \sigma(\alpha \phi), $$
hence
    $$\delta(\sigma) \phi = \alpha^{-1} \alpha \delta(\sigma)
                    \sigma\left(\sigma^{-1}(\phi)\right) 
        = \alpha^{-1}\sigma\left(\alpha\left(\sigma^{-1}(\phi)\right)\right)
              = \alpha^{-1} {^\sigma\alpha}(\phi). $$
This completes the proof of the proposition. \end{proof}
\section{Proof of the Main Theorem}
Fix $\alpha \in \GL\left({\mathcal S}(X,E)\right)$ satisfying the conclusion
of Proposition $18$. In light of $(11)$ and $(16)$, if $\sigma \in {\mathfrak H}$ 
and $g \in \Sp(V)$ then
  $$ ^\sigma\left( \alpha W_\lambda(g) \alpha^{-1} \right)
      = {^\sigma\alpha} {^{\sigma}W_\lambda(g)} {(^\sigma\alpha)^{-1}}
      =  {^\sigma\alpha} {\delta(\sigma)}^{-1} W_\lambda(g) \delta(\sigma)
            {(^\sigma\alpha)^{-1}} 
        = \alpha W_\lambda(g) \alpha^{-1}. $$
The compatability of the Galois actions $(8)$ and $(9)$ allows us
to deduce that the operators 
      $$\alpha W_\lambda(g) \alpha^{-1}, \qquad g \in \Sp(V),$$
leave
      $${\mathcal S}(X,E)^{\mathfrak H} = {\mathcal S}\left(X, E^{\mathfrak H}\right)
             = {\mathcal S}\left(X, \rn(\sqrt{p}, \sqrt{-p})\right) $$
invariant, hence provide a projective Weil representation realized
over $\rn(\sqrt{p}, \sqrt{-p})$.

\end{document}